\documentclass[11pt,a4paper]{article}

\usepackage{header}
\usepackage{abbrev}

\begin{document}

\title{Nonlinear Schwarz methods to compute geodesics on manifolds}

\author{Marco Sutti\thanks{Mathematics Division, National Center for Theoretical Sciences, National Taiwan University, Taipei, Taiwan (\email{msutti@ncts.ntu.edu.tw}).}\hspace{2mm}\orcidlink{0000-0002-8410-1372} and Tommaso Vanzan\thanks{Dipartimento di Scienze Matematiche, Politecnico di Torino, Italia (\email{tommaso.vanzan@polito.it}). T.V. is member of INdAM GNCS.}\hspace{2mm}\orcidlink{0000-0001-7554-4692}}

\date{\today}

\maketitle

\begin{abstract}
We consider the leapfrog algorithm by Noakes for computing geodesics on Riemannian manifolds. The main idea behind this algorithm is to subdivide the original endpoint geodesic problem into several local problems, for which the endpoint geodesic problem can be solved more easily by any local method (e.g., the single shooting method). The algorithm then iteratively updates a piecewise geodesic to obtain a global geodesic between the original endpoints. From a domain decomposition perspective, we show that the leapfrog algorithm can be viewed as a classical Schwarz alternating method. Thanks to this analogy, we use techniques from nonlinear preconditioning to improve the convergence properties of the method. Preliminary numerical experiments suggest that this is a promising approach.



\end{abstract}

\section{Introduction}
\label{sec:intro}
Let $\cM$ be a connected Riemannian submanifold of $\R^d$ equipped with a Riemannian metric, that is, a family of scalar products $(\langle \cdot,\cdot\rangle_p)_{p\in \cM}$ defined on the tangent spaces, $\langle\cdot,\cdot\rangle_p:\TpM\times \TpM\rightarrow \R$. In this work, we are concerned with the computation of geodesics, which are defined as the critical points of the length functional
\begin{equation}\label{eq:functional_length}
\cL(c)\coloneqq\int_0^1 \sqrt{\langle \dot{c}(t),\dot{c}(t)\rangle_{c(t)}}\;dt,
\end{equation}
$c:[0,1]\rightarrow \cM$ being a smooth curve satisfying $c(0)=p$ and $c(1)=q$.
The computation of geodesics has relevant applications in many fields since the Riemannian distance between two points on a manifold is defined as the length of the curve minimizing $\cL$. 
For some manifolds, explicit formulas exist for the minimizing geodesic and consequently for the Riemannian distance. This is the case, for instance, for the unit sphere, or the Grassmann manifold; see \protect{\cite[Theorem~8]{Wong:1967}} and \protect{\cite[\S 3.8]{AMS:2004}}. For several other manifolds, explicit formulas are not known, but geodesics are characterized as solutions of nonlinear boundary value problems, which motivates the analysis of numerical methods. As a concrete example, the boundary value problem characterizing geodesics on the Stiefel manifold is (see, e.g., \protect{\cite[(2.7)]{Edelman}})
\begin{equation}
\ddot{c}(t)+c(t)(\dot{c}(t)^\top \dot{c}(t))=0,\quad c(0)=Y_0,\; c(1)=Y_1,
\end{equation}
where $c(t)\in \mathbb{R}^{n\times k}$ is a matrix with orthonormal columns for every $t \in [0,1]$.

Numerical algorithms on manifolds often rely on two important maps in differential geometry closely related to geodesics, namely the Riemannian exponential and logarithmic maps. Given a point $p\in \cM$, the Riemannian logarithmic map $\Log_p: \cM \rightarrow \TpM$ associates to another point $q\in \cM$ the tangent vector $\dot{c}(0) \in \mathrm{T}_p{\cM}$ of the curve $c$ minimizing \eqref{eq:functional_length}. The Riemannian exponential map $ \Exp_{p}:\TpM\rightarrow \cM$ instead, given a vector $v\in \mathrm{T}_p{\cM}$, returns the point $q=c(1)$ of the curve minimizing \eqref{eq:functional_length} that satisfies $c(0)=p$ and $\dot{c}(0)=v$.  Note that computing the Riemannian logarithm is equivalent to solve the geodesic \textit{boundary value} problem (it can also be shown that the distance between $p$ and $q$ is equal to $\|\Log_p(q)\|$, $\|\cdot\|$ being the Riemannian norm on $\TpM$ \cite[Proposition~10.22]{boumal2023intromanifolds}), and therefore an analytic expression for $\Log_p$ is only rarely known. On the other hand, computing the exponential map is equivalent to solve an \textit{initial value} problem (specifying initial position and velocity), and closed formulas are often available. 

This has motivated the development of single or multiple shooting methods to approximate numerically the logarithmic map; see, e.g, \cite{Bryner,Sutti_thesis,Sutti_shot}. These methods start from an initial tangent vector $v\in \TpM$, iteratively evaluate the exponential map and correct the tangent vector, possibly within a Newton's method. These strategies however often need quite accurate initial guesses. A kind of globalization strategy was proposed in \cite{Noakes:1998}, where the author suggested to introduce intermediate points between $p$ and $q$, between which the geodesic can be computed more easily using shooting methods, and an outer relaxation procedure guarantees convergence to the global geodesic.

The contribution of this manuscript is twofold. First, in Section \ref{sec:leapfrog_Schwarz} we interpret the algorithm of \cite{Noakes:1998} as a classical alternating Schwarz method formulated in its substructured form; see, e.g., \cite{Ciaramella_Vanzan,Ciaramella_Vanzan2}. Second, this new point of view allows us to apply many strategies developed by the domain decomposition community to speed up the convergence. Here, we start considering nonlinear preconditioning in Section \ref{sec:RASPEN}, and we present promising numerical results on the sphere and on the Stiefel manifold in Section \ref{sec:num_ex}. 

\section{The leapfrog method as a Schwarz method}
\label{sec:leapfrog_Schwarz}
\begin{figure}[htbp]
   \centering
      \includegraphics[width=0.27\columnwidth]{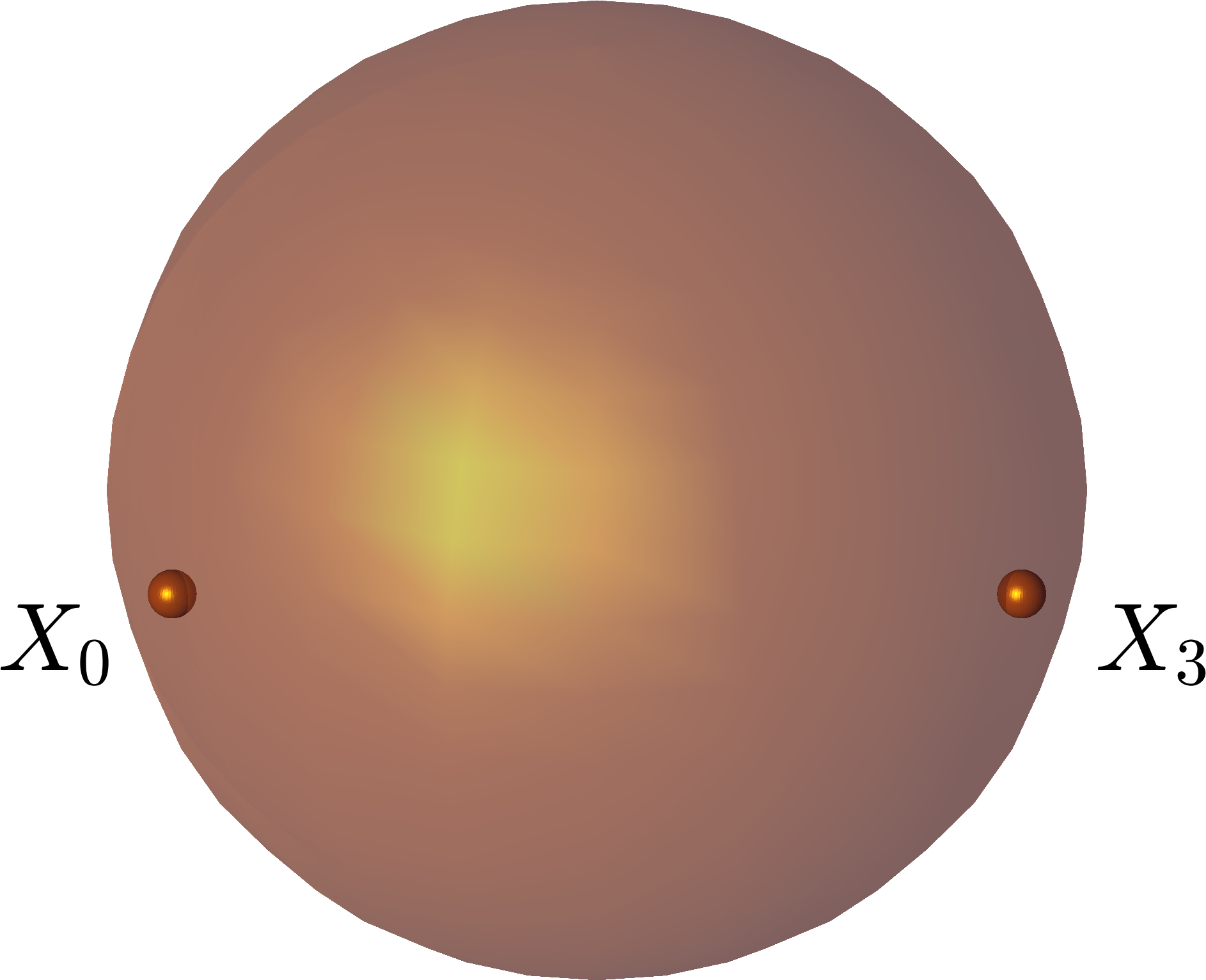} \hspace{5mm}
      \includegraphics[width=0.27\columnwidth]{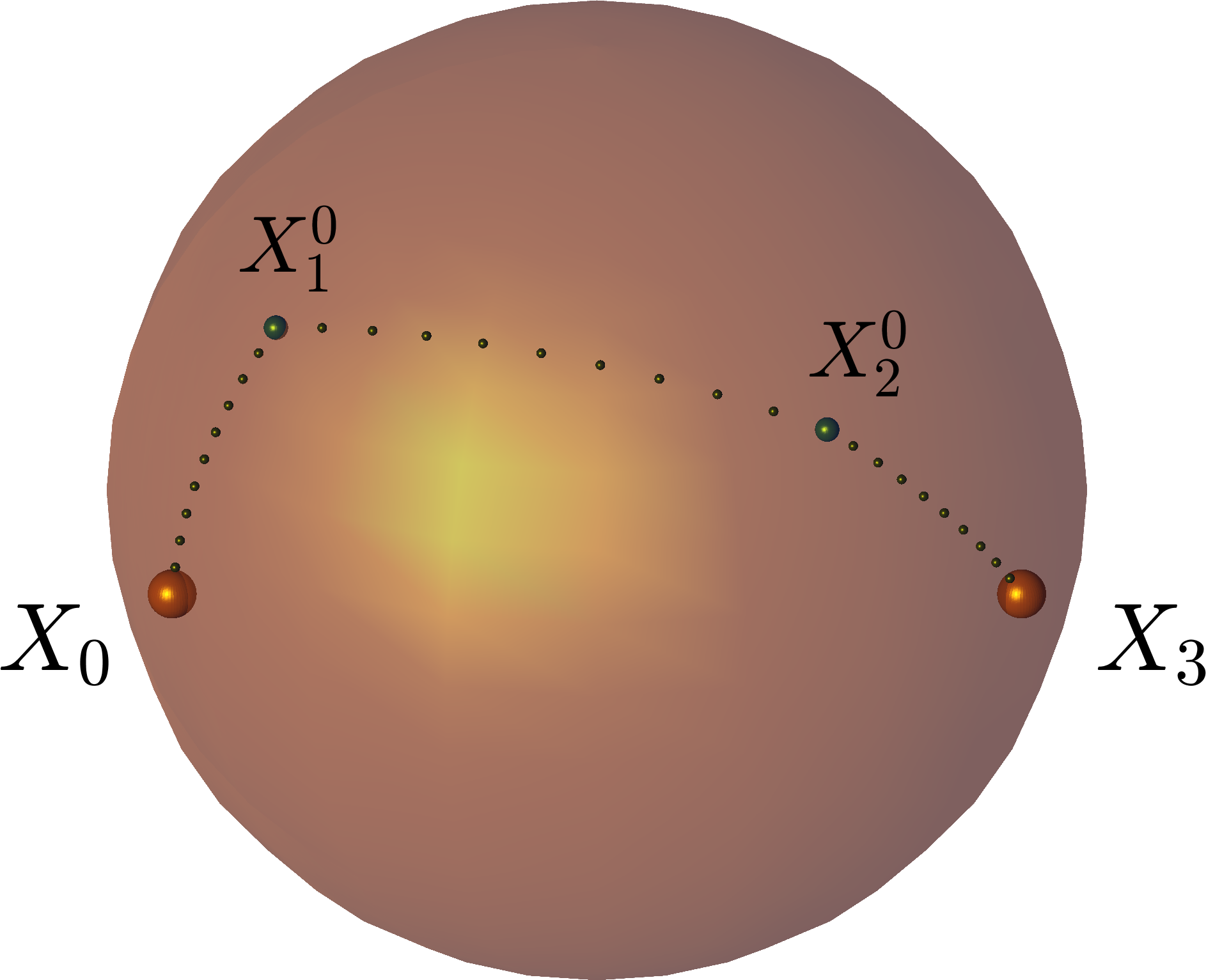} \hspace{5mm}
      \includegraphics[width=0.27\columnwidth]{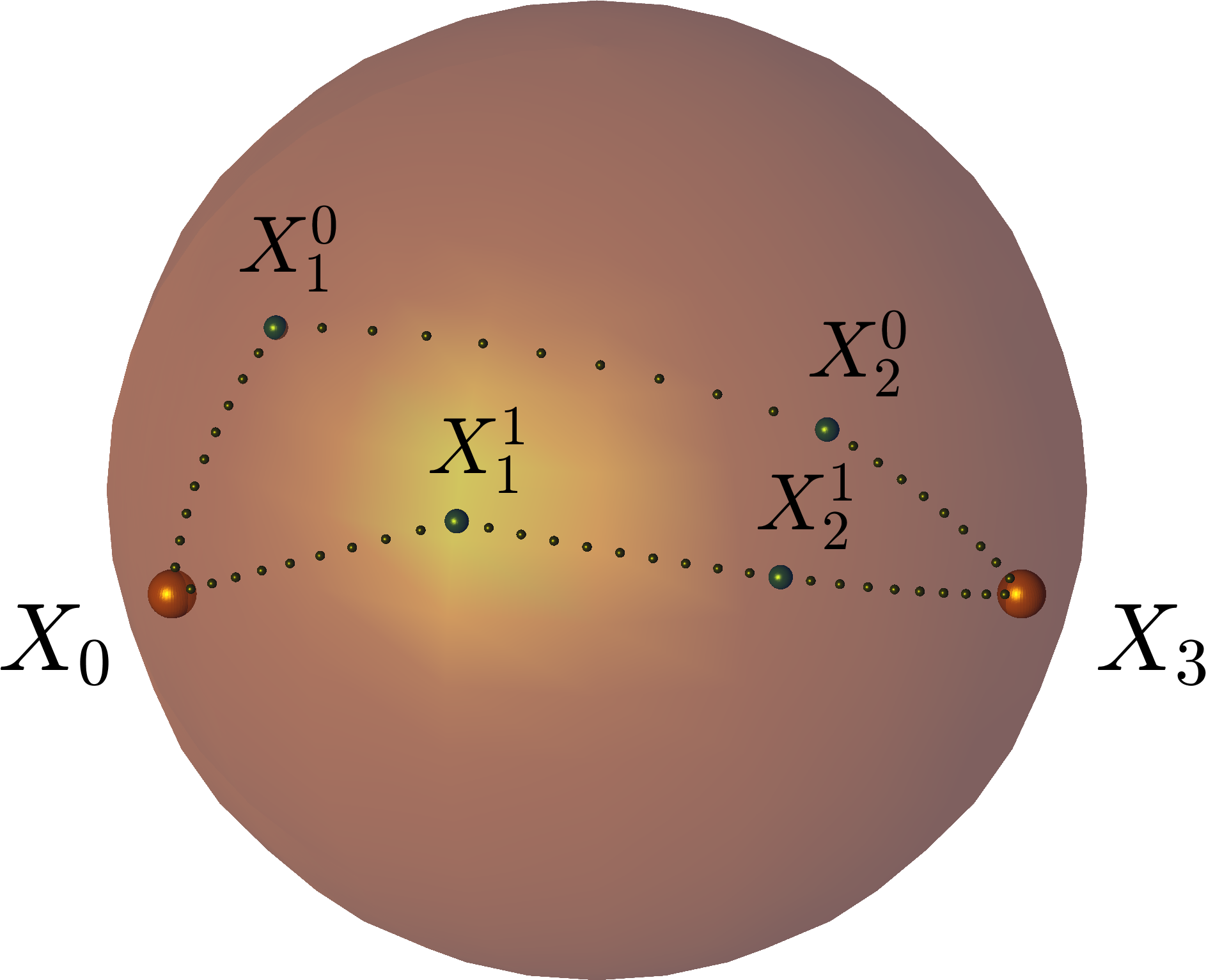}\\\vspace{7mm}
      \includegraphics[width=0.27\columnwidth]{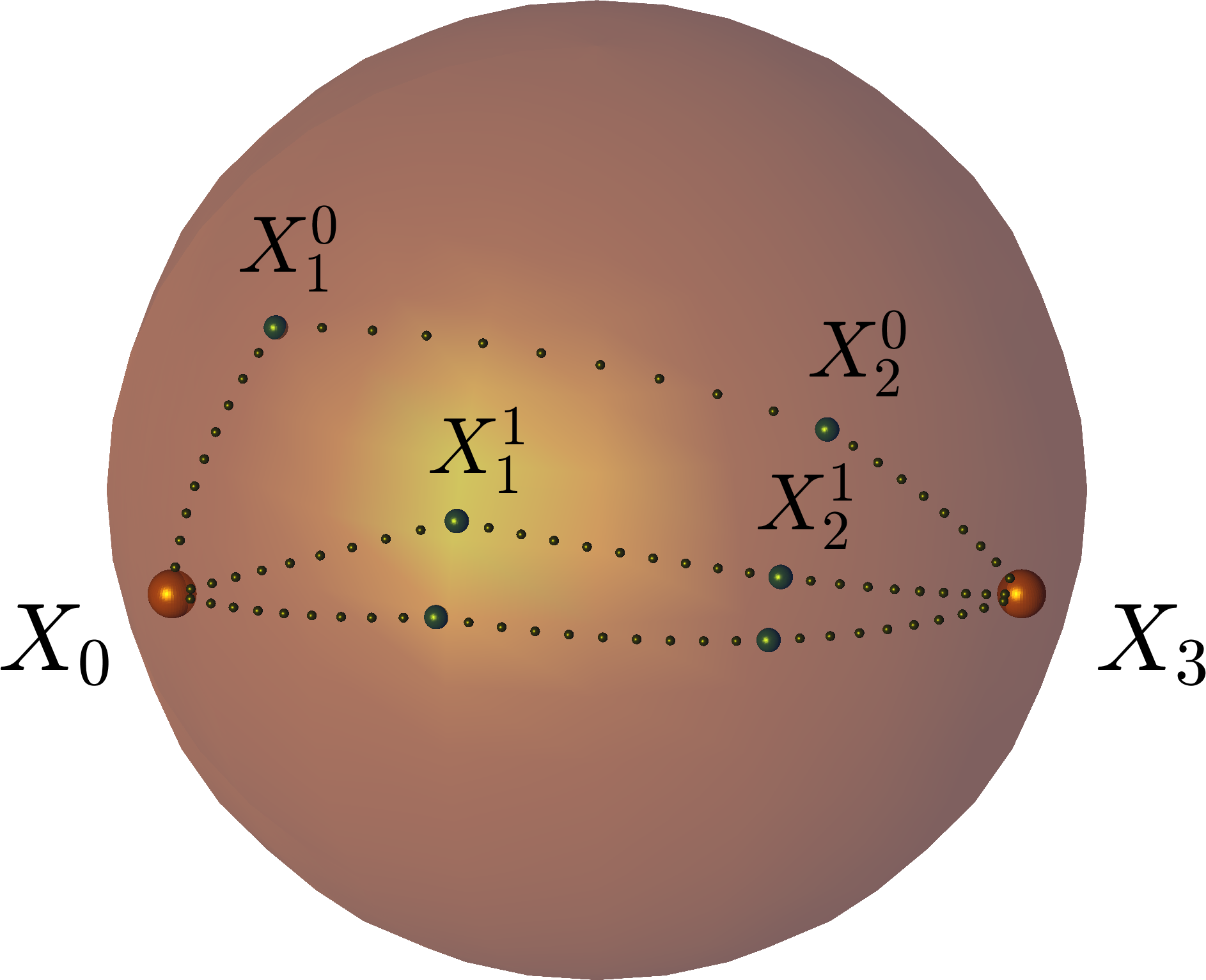} \hspace{5mm}
      \includegraphics[width=0.27\columnwidth]{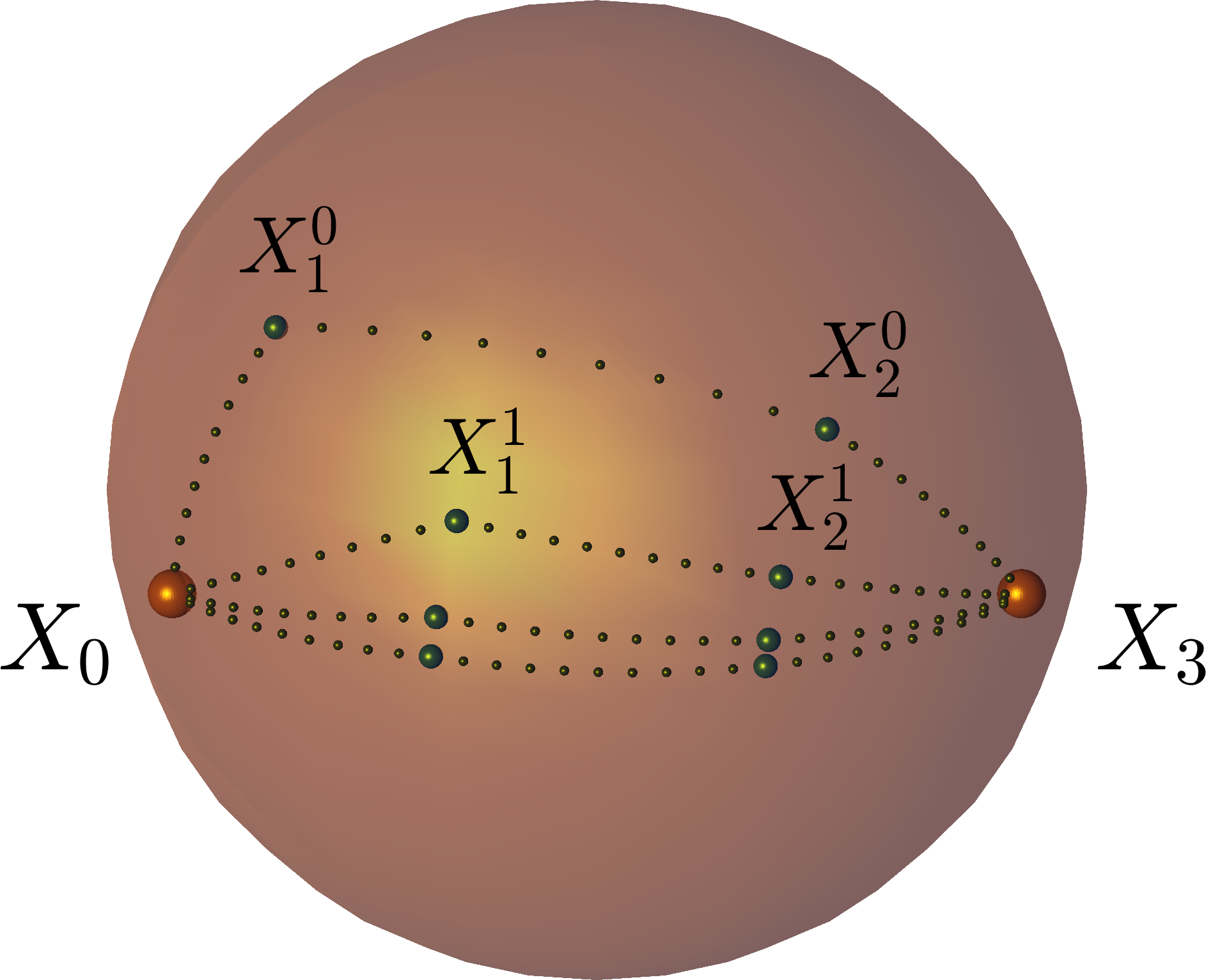} \hspace{5mm}
      \includegraphics[width=0.27\columnwidth]{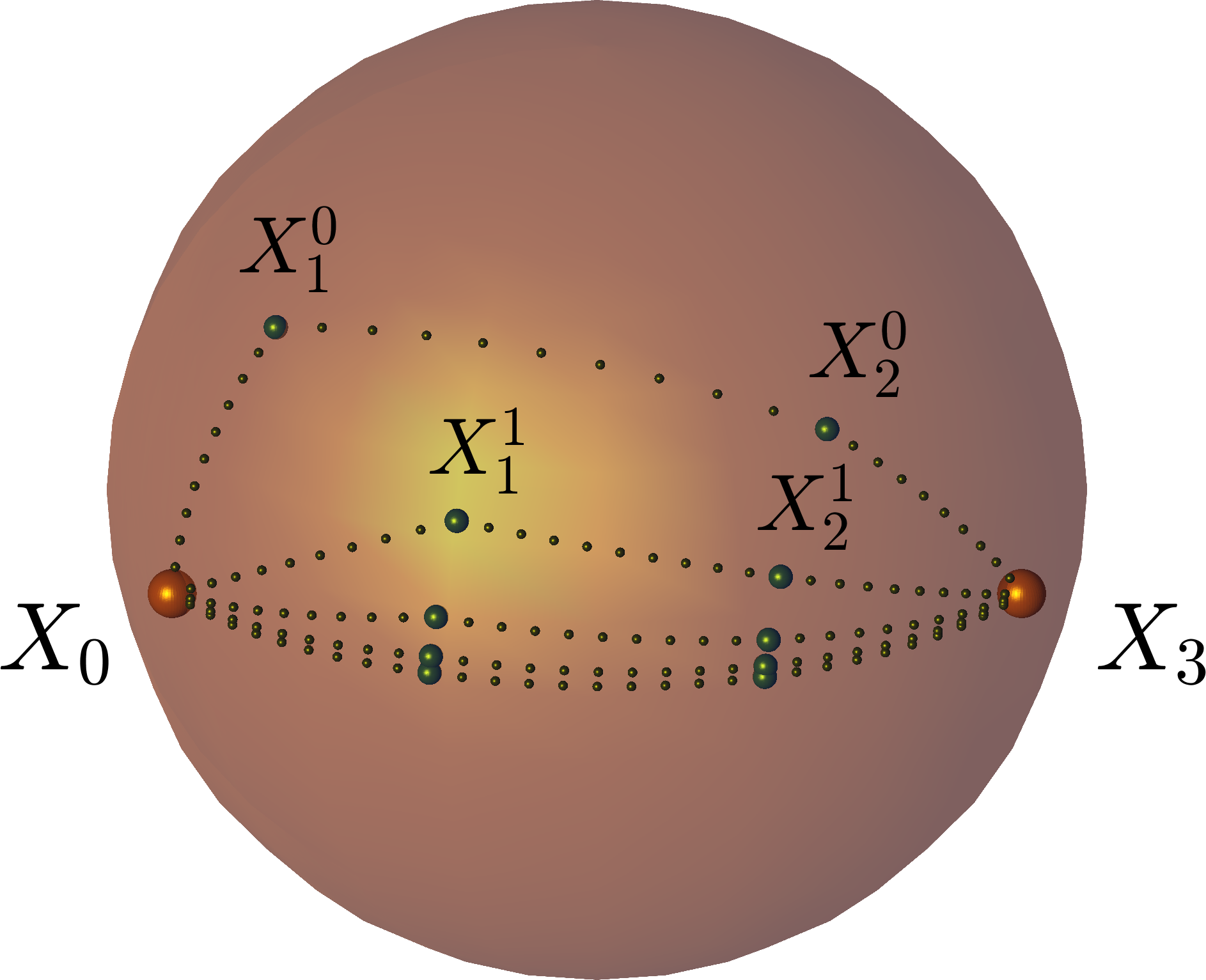}
      \caption{Illustration of the first iterates of the leapfrog algorithm with $m=4$.}
   \label{fig:leapfrog_procedure}
\end{figure}

The algorithm proposed by Noakes in \cite{Noakes:1998} is a classical instance of the \emph{think globally, act locally} paradigm that typically characterizes domain decomposition methods. As previously mentioned, given two points $p$ and $q$ on $\cM$, the direct computation of a geodesic using a shooting method may easily fail if the initial guess is not sufficiently accurate. To overcome this limitation, Noakes suggested to consider the $m$-tuple $\widetilde{\Xb}=(X_0,X_1,\dots,X_{m-2},X_{m-1})^\top\in \cM^{m}$, $m\in \mathbb{N}$, with $X_0=p$ and $X_{m-1}=q$, and the associated \textit{piecewise geodesic} obtained by joining the geodesics that link the couple $(X_i,X_{i+1})$ for every $i\in\left\{0,\dots,m-2\right\}$\footnote{The motivating idea is that if the points $(X_i,X_{i+1})$ are close enough to each other, a shooting algorithm may exhibit a faster and more robust convergence with respect to the choice of the initial tangent vector.}.
The desired geodesic connecting $X_0$ and $X_{m-1}$ is then obtained within an iterative algoritm: 
let $\widetilde{\Xb}^k=(X_0,X^k_1,\dots,X^k_{m-2},X_{m-1})^\top$ be the current iterate at iteration $k$, then a geodesic is computed between the points $(X_0,X^{k}_{2})$, and $X^{k+1}_{1}$ is set equal to the midpoint of that geodesic; see Figure \ref{fig:leapfrog_procedure} for a graphical description.
Then, in a Gauss--Seidel fashion, $X^{k+1}_{2}$ is set equal to the midpoint of the geodesic connecting $X^{k+1}_1$ and $X^{k}_{3}$, and this process is repeated until $X^{k+1}_{m-2}$ is computed.
Noakes showed that the algorithm generates a sequence of piecewise geodesics of decreasing length and, by a compactness argument, that there exists at least a convergent subsequence such that all the points of the limit curve lie on the global geodesic connecting $X_0$ and $X_{m-1}$. More recently, \cite{Sutti_thesis} analyzed the convergence of the algorithm interpreted as a nonlinear Gauss--Seidel minimization algorithm.
The leapfrog algorithm of \cite{Noakes:1998} is summarized by the pseudocode in Algorithm \ref{Alg:Leapfrog}. A parallel Jacobi version can be straightforwardly derived, although it less common as the decreasing length property of the piecewise geodesic may not hold.

\begin{center}
\begin{minipage}{0.90\textwidth}%
\begin{algorithm}[H]
   \SetAlgoLined
   \KwData{Two points $X_{0}$, $X_{m-1}$, an integer $m$.}
   \KwResult{Geodesic connecting $X_{0}$ and $X_{m-1}$.}
   Compute the initial guess for the $m-2$ intermediate points\;
   $k = 0$\;
   \While{a stopping criterion is not met}{
   \For{$i=1 \colon m-2$}{
	  Compute the geodesic between $X_{i-1}^{k+1}$ and $X_{i+1}^k$ \;
	  Set $X_{i}^{k+1}$ equal to the midpoint\;    
      }
      Update $k \leftarrow k + 1$\;
   }
\caption{The leapfrog algorithm.}\label{Alg:Leapfrog}
\end{algorithm}
\end{minipage}
\end{center}

To formalize this algorithm within a Schwarz framework, let us introduce a time discretization $t_0<t_1<\dots <t_{m-1}$ of the interval $[0,1]$ with $t_0=0$, $t_{m-1}=1$ and $t_{i}=\frac{t_{i-1}+t_{i+1}}{2}$. The points $X_i$ are then thought to be guessed values of the sought geodesic $c(t)$ at times $t_i$. Next, we introduce the midpoint map $\mathfrak{M}\colon \cM \times \cM \to \cM $, defined by 
\[\Mid(x,y) = \Exp_{x} \Big( \tfrac{1}{2}\,\Log_{x}(y) \Big),\]
which takes two points $x,y\in \cM$, computes the tangent vector of the geodesic joining $x$ and $y$, and then moves through the exponential map half way along the geodesic. 
We then make the key observation\footnote{The first remark in this direction was made by Martin J. Gander during the PhD defence of the first author.} that the operator $\Mid$ precisely acts as a subdomain solver: it receives two boundary values, say $X_{i-1}^{k+1}$ and $X_{i+1}^k$, solves the nonlinear boundary value problem within the action of the logarithmic map, and then it returns the value of the subdomain solution (curve) at time $t_i$, i.e., precisely in the middle of the interval $[t_{i-1},t_{i+1}]$.
The algorithm can then be formulated as a substructured Schwarz method \cite{Ciaramella_Vanzan,Ciaramella_Vanzan2}, that is, a Schwarz method formulated as an iteration directly on the trace values: given an iterate $\Xb^k=(X_1^k,\dots,X_{m-2}^k)^\top$, one iteration of (the Jacobi variant of) the leapfrog algorithm can be written compactly as the Schwarz iteration
\begin{equation}\label{eq:Schwarz_iteration}
\Xb^{k+1}=\begin{pmatrix}
X_1^{k+1}\\[4pt]
X_2^{k+1}\\[4pt]
\vdots\\[4pt]
X_{m-3}^{k+1}\\[4pt]
X_{m-2}^{k+1}
\end{pmatrix}=\begin{pmatrix}
\Mid(X_0,X_2^{k})\\[4pt]
\Mid(X_1^k,X_3^{k})\\[4pt]
\vdots\\[4pt]
\Mid(X_{m-4}^k,X_{m-2}^{k})\\[4pt]
\Mid(X_{m-3}^k,X_{m-1}^{k})
\end{pmatrix}
=\G\left(\Xb^{k}\right).
\end{equation}

\section{A preconditioned nonlinear Schwarz method}
\label{sec:RASPEN}
In the field of nonlinear preconditioning, a key idea is to apply Newton's method to the fixed point equation associated to the Schwarz iteration, see, e.g., \cite{ASPIN,Dolean,Chaouqui}. Here, the fixed point equation associated to \eqref{eq:Schwarz_iteration} can be described by the nonlinear equation $\F(\Xb)=0$, where $\F:\cM^{m-2}\rightarrow \R^{d(m-2)}$ is defined as $\F(\Xb)\coloneqq \Xb-\G\left(\Xb\right)$,
and whose Jacobian $J_{\F(\Xb)}$, due to the particular structure of $\G$, has the block tridiagonal form
\begin{small}
\begin{equation}\label{eq:Jacobian}
   J_{\F(\Xb)} = I - J_{\G(\Xb)}=
      \begin{bmatrix}
         I  & -g^{0,2}_2  &  &  &  \\[12pt]
         -g^{1,3}_1  & I  & -g^{1,3}_3 &  & \\[12pt]
         & \ddots & \ddots & \ddots &  \\[12pt]
         & & -g^{m-3,m-1}_{m-3} & I & -g^{m-3,m-1}_{m-1}  \\[12pt]
         & & & -g^{m-2,m}_{m} & I  \\
      \end{bmatrix},
\end{equation}
\end{small}
where $g^{i,j}_{k} \coloneqq \partial_{X^{k}} \mathfrak{M}(X_{i},X_{j})$. 
In our setting though, due to the nonlinear nature of $\cM$, it is not obvious how to perform a Newton method to compute a zero of $\F$ since in general $\Xb^k+\deltab \notin \cM^{m-2}$, $\deltab$ being the solution of the Newton system. To overcome this, we consider a modified Newton method, described in Algorithm~\ref{Alg:modifiedNewton}, which differs in two additional steps. First, the Newton direction $\deltab$, which could be any vector in the embedding space $\R^{d(m-2)}$, is projected onto the tangent space $\mathrm{T}_{\Xb^k}\cM^{m-2}$ to $\cM^{m-2}$ at $\Xb^k$, obtaining a tangent vector $\xib$ \cite[\S 3]{boumal2023intromanifolds}. Second, using a retraction mapping $\cR_{\Xb^k}$ at $\Xb^k$, we move along the direction identified by $\xib$ while remaining constrained to the manifold $ \cM^{m-2} $; see, e.g., \cite[\S 3.6]{boumal2023intromanifolds}.
\begin{center}
\begin{minipage}{0.90\textwidth}%
\begin{algorithm}[H]
   \SetAlgoLined
   \KwData{Initial guess $\Xb^0$.}
   \KwResult{$\Xb$ such that $\F(\Xb)=0$.}
   Set $k=0$\;
   \While{a stopping criterion is not met}{
    Solve $\delta=-J_{\F(\Xb^k)}^{-1}\F(\Xb^k)$\;
    Set $\xib:=\cP_{\mathrm{T}_{\Xb^k}\cM^{m-2}}(\deltab)$\;
    Set $\Xb^{k+1}=\cR_{\Xb^k}\left(\Xb^{k+1}+\xib\right)$\;
    Update $k=k+1$.
      }
\caption{The modified Newton algorithm.}\label{Alg:modifiedNewton}
\end{algorithm}
\end{minipage}
\end{center}
Incidentally, we mention that an alternative method would be to consider the Riemannian Jacobian $J^{R}_{\F(\Xb)} \colon \mathrm{T}_{\Xb}\cM^{m-2}\rightarrow \R^{d(m-2)}$, and then to look for a tangent vector $\xib$ such that $J^{R}_{\F(\Xb)}\xib\approx \F(\Xb)$ (e.g., in a least square sense). One could even minimize the norm of the residual of $\F(\Xb)$ in a Gauss--Newton framework, closer to the spirit of \cite{Noakes22}.
However, in this manuscript we limit ourself to Algorithm~\ref{Alg:modifiedNewton}.

\section{Numerical examples}
\label{sec:num_ex}
In this section, we compare the behavior of leapfrog and of the preconditioned Schwarz method (hereafter denoted preconditioned leapfrog) for computing geodesics between two points on $\cM$.
We first consider the case of the unit sphere $ \cS^{d-1}\coloneqq \lbrace \x \in \R^{d} \colon \x\tr \x = 1 \rbrace$ which is a Riemannian submanifold of $ \mathbb{R}^{d} $. Even though for $\cS^{d-1}$ there is an explicit formula for computing the logarithmic map (and thus the Riemannian distance between two points), we consider these experiments to be preliminary before moving to manifolds for which closed-form formulas are not available.
We briefly recall the geometry of $ \cS^{d-1}$. The tangent space to $ \cS^{d-1} $ at $ \x $ is 
$\mathrm{T}_{\x} \cS^{d-1} \coloneqq \lbrace \z \in \R^{d} \colon \x\tr \z = 0 \rbrace $, equipped with an inner product inherited from the embedding space $ \R^{d} $, i.e., for $ \boldsymbol{\xi}, \boldsymbol{\eta} \in \mathrm{T}_{\x} \cS^{d-1} $, $ \langle \boldsymbol{\xi}, \boldsymbol{\eta} \rangle_{\x} := \boldsymbol{\xi}\tr \boldsymbol{\eta} $.
The projector $\cP_{\x}$ onto $ \mathrm{T}_{\x} \cS^{d-1} $ is
$ \cP_{\x} \boldsymbol{\xi} = (I_{d}-\x\x\tr)\,\boldsymbol{\xi} $, while the retraction $\mathcal{R}_{\x}$ is defined as $\mathcal{R}_{\x}(\boldsymbol{\xi}):=\frac{\x+\boldsymbol{\xi}}{\|\x+\boldsymbol{\xi}\|}$.
The Riemannian exponential at $\x$ takes a tangent vector $\dt{\x}\in \mathrm{T}_{\x}\cS^{d-1}$ and returns a point $ y \in \cM $ such that
\begin{equation}\label{eq:sphere_exponential}
   \y = \Exp_{\x}(\dt{\x}) = \x \cos(\| \dt{\x} \|) + \frac{\dt{\x}}{\| \dt{\x} \|} \sin(\| \dt{\x} \|  ),
\end{equation}
while the corresponding Riemannian logarithm is
\begin{equation}\label{eq:sphere_logarithm}
   \Log_{\x}(\y) = \dt{\x} = \arccos(\x\tr \y) \frac{\cP_{\x}\,\y}{\| \cP_{\x}\,\y \|}.
\end{equation}
Thanks to these analytical expressions, in this test case we are able to compute by hand the terms $g_{k}^{i,j}$, and thus we fully assemble the Jacobian $J_{\F(\Xb)}$ in \eqref{eq:Jacobian}.

\begin{figure}[]
   \centering
      \includegraphics[height=3cm,width=\columnwidth]{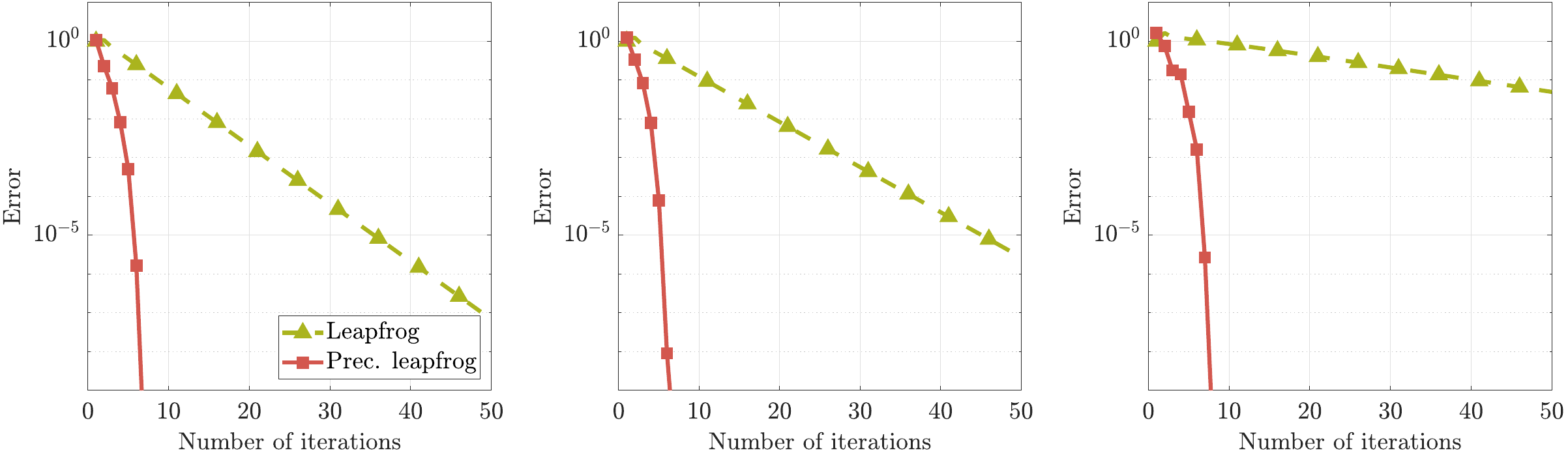}
      \caption{Unit sphere with $ d = 100 $, $ m = 4 $, and increasing distance $ \mathrm{dist}(X,Y) = 0.1\pi, 0.5\pi, 0.9\pi $.}
   \label{fig:Sphere_n100_m4}
\end{figure}

Figure~\ref{fig:Sphere_n100_m4} shows the convergence behavior of leapfrog and its preconditioned version for $\cS^{d-1}$ with $ d = 100 $, $ m = 4 $, and increasing distance between the endpoints $ \mathrm{dist}(X,Y) = 0.1 \pi, 0.5 \pi$ and  $0.9\pi $. It is evident that the leapfrog method gets slower as the distance between the endpoints increases, while preconditioned leapfrog exhibits a quadratic convergence behavior which is basically independent of the distance between the endpoints. In all the cases, leapfrog clearly requires many more iterations to reach the same accuracy level as its preconditioned version. In addition, notice that leapfrog is sequential, while its preconditioned version, being based on the Jacobi Schwarz iteration \eqref{eq:Schwarz_iteration}, may take further advantage of parallelization.

\begin{figure}[]
   \centering
      \includegraphics[height=3cm,width=\columnwidth]{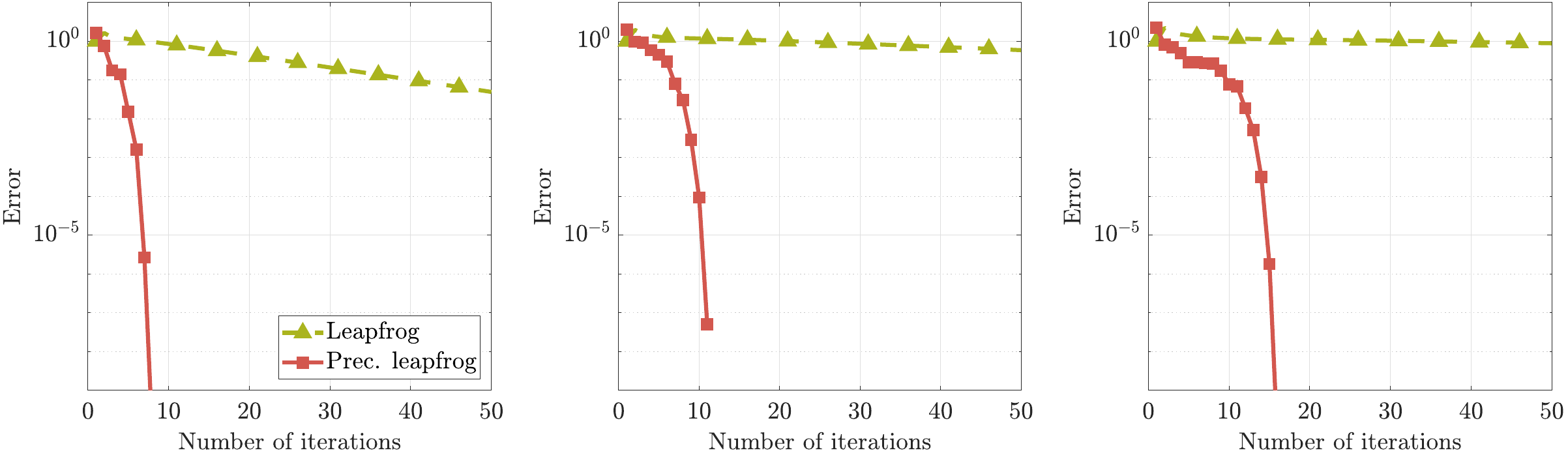}
      \caption{Unit sphere with $ d = 100 $, $ \mathrm{dist}(X,Y) = 0.9\pi $, and varying number of points $ m = 4, 7, 10 $.}
   \label{fig:Sphere_n100_dist0.9pi}
\end{figure}

Figure~\ref{fig:Sphere_n100_dist0.9pi} shows the convergence behavior of leapfrog and its preconditioned version for a fixed distance between the endpoints $\mathrm{dist}(X,Y)=0.9\pi$ and increasing values of $m$ (i.e., of subdomains). In this setting, the convergence of leapfrog significantly deteriorates, while preconditioned leapfrog exhibits a much weaker dependence on $m$.

We next consider the Stiefel manifold $ \mathrm{St}(n,p)\coloneqq\left\{Q\in \R^{n\times p}:\; Q^\top Q=I_p\right\} $ of orthonormal $n$-by-$p$ matrices. 
Due to space limitation, we do not discuss the geometry of the Stiefel manifold here, and we refer the reader to, e.g., \cite{Edelman,boumal2023intromanifolds,Sutti_shot}. Since for the Stiefel manifold there is no closed expression for the logarithmic map, the midpoint map $\mathfrak{M}$ is evaluated with a numerical routine (e.g., a single shooting method or Bryner's method \cite{Bryner}), and we solve the Jacobian system in a matrix-free fashion by computing the action of the Jacobian onto a vector with a finite difference approximation.

\begin{figure}[]
   \centering
      \includegraphics[height=3cm,width=\columnwidth]{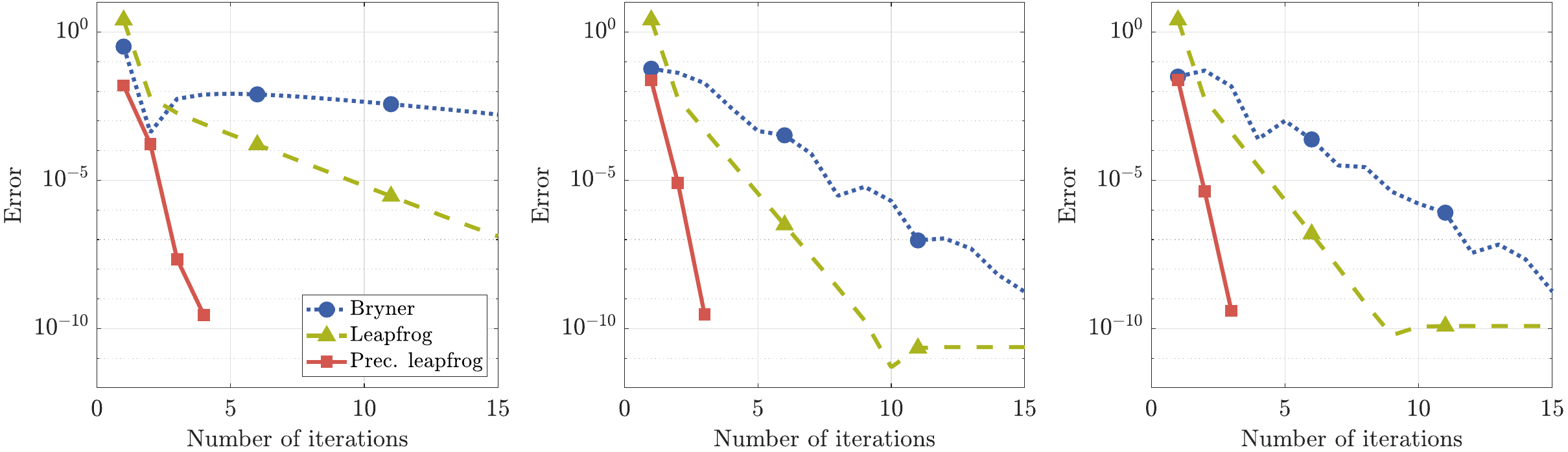}
      \caption{Stiefel manifold with $ n = 100 $, $ m = 4 $, $ \mathrm{dist}(X,Y) = 0.8\pi $, and increasing $ p = 2, 12, 22 $.}
   \label{fig:Comparison_Methods_St_n100_varying_p_distXY_0.8pi}
\end{figure}

Figure~\ref{fig:Comparison_Methods_St_n100_varying_p_distXY_0.8pi} compares the convergence behavior of three numerical methods, namely leapfrog, Bryner's method, and preconditioned leapfrog for a Stiefel manifold with $ n = 100 $, $ m = 4 $, distance between the endpoints of $ \mathrm{dist}(X,Y) = 0.8\pi $, and for increasing values of $p$, namely $ p = 2, 12, 22 $.
Preconditioned leapfrog seems to be independent of $p$, while for leapfrog and Bryner's methods the problem appears to be more difficult when $ p \ll n $. This is in agreement with similar observations made in \protect{\cite[Table 5.1]{Zimmermann:2017}}, \protect{\cite[\S 5.2]{Nguyen:2022}}, and \protect{\cite[\S 5.1]{Sutti_shot}}.

\section{Conclusions}
We interpreted a classical algorithm to compute geodesic on manifolds as a classical Schwarz method. This allowed us to consider nonlinear preconditioning strategies which effectively speed up the convergence. 
There are several research directions that are now opening. First, in concrete problems, like for the Stiefel manifold, the midpoint map is not available in closed form, and thus the computation of the action of the Jacobian onto a vector is delicate.
Future efforts will focus on the use of implicit differentiation (as suggested in \cite{Noakes22}) to evaluate the Jacobian action based on the derivative of the exponential map, or on quasi-Newton approximations. 
Second, as mentioned, there are several variants of Newton's method, depending on how much geometry one wishes to enforce, and it is not clear at the moment if one is to be preferred. Third, coarse spaces may be developed to improve the scalability of both the standard and preconditioned leapfrog with respect to $m$.
Finally, applications to other manifolds deserve to be further investigated.





\end{document}